\documentclass{article} 

\usepackage{amsfonts,amsmath,amssymb}
\usepackage[all]{xy} 
\usepackage{epsfig, latexsym}

\usepackage[french]{babel}
\usepackage[latin1]{inputenc}

\newenvironment{preuve}[1][]
{\vskip 2mm\noindent {\bf D\'emonstration #1:  }}{$\Box$ \vskip 2mm}
\newenvironment{theor}
{\vskip 2mm\noindent {\bf Th\'eor\`eme}}{\vskip 2mm}
\newcommand{\db}{\bar\partial}
\newcommand{\Nn}{\mathbb{N}}
\newcommand{\Rr}{\mathbb{R}}

\newcommand{\Zz}{\mathbb{Z}}
\newcommand{\Cc}{\mathbb{C}}
\newcommand{\ep}{\epsilon}
\newcommand{\om}{\omega}
\newcommand{\lam}{\lambda}

\newtheorem{lemma}{Lemme}

\newtheorem{prop}{Proposition}

\title{Convexit\'e rationnelle des sous-vari\'et\'es immerg\'ees lagrangiennes}

\begin{document}

\author {\normalsize{Damien Gayet }}

\maketitle
\centerline{\textbf{Abstract}}
We prove that a  compact, immersed, submanifold of $\Cc^n$, lagrangian for
a K\"{a}hler form, is rationally convex, generalizing a theorem
of Duval and Sibony for embedded submanifolds.\\

\centerline{\textbf{R\'esum\'e}}
Nous d\'emontrons qu'une sous-vari\'et\'e compacte, immerg\'ee dans $\Cc^n$ et 
lagrangienne pour une forme de K\"{a}hler, est rationnellement
convexe, g\'en\'eralisant ainsi un th\'eor\`eme de Duval et Sibony
pour des sous-vari\'et\'es plong\'ees.\\

\textsc{Mots-clefs}: Convexit\'e rationnelle, Vari\'et\'e lagrangienne immerg\'ee,
Estim\'ees $L^2$ de H\"{o}rmander.

\textsc{Code mati\`ere AMS}: 32E20-32F20

 \section{Introduction}
Soit $X$ un compact dans $\Cc^n$.  On consid\`ere l'enveloppe rationnelle
$r(X)$ de $X$. Un point $x$ est dans $r(X)$ si et seulement
si toute hypersurface alg\'ebrique passant
par $x$ rencontre $X$. Le compact $X$ est \textit{rationnellement convexe}
si $r(X)= X$, i.e si son compl\'ementaire est r\'eunion d'hypersurfaces \'evitant $X$.
 Un des int\'er\^ets de cette notion r\'eside
dans la g\'en\'eralisation du th\'eor\`eme de Runge:
\begin{theor}\textit{
(Oka, Weil) Toute fonction holomorphe au voisinage
d'un compact rationnellement convexe peut \^etre approch\'ee
uniform\'ement sur le compact par des fractions rationnelles.}
                                
        \end{theor}
  Une obstruction classique \`a la convexit\'e rationnelle repose sur
le fait suivant:
\textit{Une surface de Riemann $\Sigma$ compacte \`a bord, telle que
$\partial \Sigma$ borde une surface
$V$ dans $X$, est contenue dans $r(X)$}.

 En effet, soit $C$ une hypersurface passant par un point de $\Sigma$.
 Les intersections de $C$ avec $\Sigma$ sont toujours positives,
 donc l'intersection homologique des deux ensembles analytiques
 est non nulle. La surface
  form\'ee par la r\'eunion de $\Sigma$ et
 de $V$ \'etant ferm\'ee,  l'intersection
  homologique de $C$ avec $\Sigma \cup V$ est
 nulle. Par cons\'equent $C$ rencontre  $X$. 
 
 Une forme $\om$  de bidegr\'e $(1,1)$
dans $\Cc^n$ est  de K\"{a}hler si elle est ferm\'ee et strictement positive.
Une sous-vari\'et\'e r\'eelle $S$ de $\Cc^n$
est dite isotrope (lagrangienne si $\dim_\Rr S =n$) pour $\om$
si la forme s'annule tangentiellement \`a $S$. 
En d'autres termes, $j^*\om=0$, o\`u $j: S \to \Cc^n$ est l'inclusion de $S$ dans
$\Cc^n$.
La sous-vari\'et\'e $S$ est alors n\'ecessairement totalement r\'eelle, i.e 
jamais tangente \`a une droite complexe.
Si $S$ est isotrope pour $\omega$,
 il n'existe pas de surface de Riemann du type pr\'ec\'edent. En effet,
 puisque $\om$ est exacte,
  $\int_V \om = - \int_\Sigma \om$, o\`u la premi\`ere int\'egrale
  est nulle, et la seconde strictement n\'egative. \\

En 1995, J. Duval et N. Sibony ([Du,Si], voir aussi [Du] dans le 
cas des surfaces) prouvent
le r\'esultat suivant:
\begin{theor} \textit{
Une sous-vari\'et\'e S compacte et lisse, plong\'ee dans $\Cc^n$, isotrope pour
une forme de K\"{a}hler, est rationnellement convexe.
                                }
 \end{theor}
 
 La r\'eciproque, facile, est vraie: une sous-vari\'et\'e totalement r\'eelle
rationnellement convexe
 est isotrope pour une certaine forme
de K\"{a}hler.\\

Il est naturel de s'int\'eresser aux sous-vari\'et\'es
lagrangiennes compactes seulement immerg\'ees dans $\Cc^n$, car
les plongements lagrangiens sont relativement rares.  Par
exemple, dans son article fondateur [Gr] de 1985,  Gromov  prouve que dans $\Cc^n$,
il n'existe pas de sous-vari\'et\'e ferm\'ee exacte
pour la forme standard $\omega_0$ (i.e $j^*d^c |z|^2$ est
non seulement ferm\'ee sur $S$ mais aussi exacte). En particulier, il n'existe pas
de sous-vari\'et\'e ferm\'ee simplement connexe, et
donc de sph\`ere $S^n$ pour $n\ge 2$, lagrangienne plong\'ee
dans $(\Cc^n, \om_0)$.  Cependant les sph\`eres lagrangiennes
\textit{immerg\'ees}  existent.  L'exemple le plus
simple est l'immersion de Whitney:

\begin{eqnarray*}
j:  & S^n \to & \Cc^n \\ & (x,y) \mapsto & j(x,y) = (1+iy)x\\
\end{eqnarray*}
o\`u $S^n$ est la sph\`ere unit\'e de $\Rr^n \times \Rr$.
La  sph\`ere immerg\'ee $j(S^n)$ pr\'esente un point double en $0$, et est
lagrangienne pour la forme standard.

Nous traiterons ici le cas des sous-vari\'et\'es lagrangiennes compactes immerg\'ees
\textit{g\'en\'eriques}, qui  poss\`edent alors un nombre fini de points
doubles ordinaires (transverses). Nous obtenons ainsi le:
\begin{theor} \textit{
 Une sous-vari\'et\'e de dimension $n$ compacte, lisse, immerg\'ee  dans $\Cc^ n $,
        lagrangienne pour une forme de K\"{a}hler $ \om $,
        et poss\'edant un nombre fini de points doubles transverses, est
        rationnellement convexe.
                                }
\end{theor}

\textit{Remarque: } En fait, le r\'esultat est valable 
dans le cas de points multiples avec  un nombre fini de
branches transverses.

\subsection{R\'esum\'e de la d\'emonstration du cas plong\'e}

Soit $S$ une sous-vari\'et\'e compacte 
lagrangienne pour une forme de K\"{a}hler
$\om$ dans $\Cc^n$.  La d\'emonstration pr\'esent\'ee dans [Du, Si] pour
prouver que $S$ est rationnellement convexe s'effectue en deux grandes
\'etapes. La premi\`ere  consiste \`a prouver  par les estim\'ees $L^2$ de
H\"{o}rmander qu'on peut remplacer, dans la d\'efinition de la convexit\'e
rationnelle, "hypersurface alg\'ebrique" par "support d'une $(1,1)$-forme
positive ferm\'ee". Cela r\'esulte de l'approximation, non seulement au sens
des courants, mais aussi au sens de la distance de Hausdorff des supports,
d'une  $(1,1)$-forme positive ferm\'ee par des hypersurfaces.

La deuxi\`eme
\'etape consiste, pour tout
$x \notin S$, \`a exhiber une $(1,1)$-forme $\om_x$ positive ferm\'ee, de
support contenant $x$ et ne rencontrant pas $S$. Ces formes s'obtiennent
par exemple en  construisant une
famille de fonctions  $\phi_\ep$ plurisousharmoniques sur $\Cc^n$,
pluriharmoniques au voisinage de $S$ et strictement plurisousharmoniques
hors des  voisinages tubulaires
 $S_\ep =\{d(.,S)\le \ep\}$. 
Ces fonctions sont des d\'eformations
d'un potentiel $\phi $ de la forme $\om$, i.e d'une fonction $\phi$
satisfaisant $dd^c\phi=i d(\db-\partial)\phi = \om$. Dans ce but, on cherche
\`a tirer profit de l'annulation tangentielle de $dd^c\phi$ le long de $S$.
On veut d'abord cr\'eer une fonction proche de $\phi$, dont le $dd^c$ est
nul (dans toutes les directions cette fois) \`a un ordre assez grand sur
$S$. On voudrait ensuite propager l'annulation du $dd^c$ \`a un
petit voisinage de $S$. 
Pour cela, en occultant les probl\`emes dus \`a
l'homologie de $S$, la d\'emarche de [Du,Si]  consiste \`a
\'etendre $\phi+i\psi$ hors de $S$ en
une fonction $h$ $\db$-plate sur $S$, o\`u $\psi$ est
 une primitive de
$j^*d^c \phi $. 
La diff\'erence $\phi - \Re h$ cro\^{\i}t alors automatiquement
en $d^2(.,S)$, car $\phi$ est strictement plurisousharmonique. 
On r\'esoud ensuite $\bar{\partial} h = \bar{\partial} u_\ep$
sur le voisinage  tubulaire $S_\ep$ (qui est pseudoconvexe car $S$ est
totalement r\'eelle), avec des estim\'ees $L^2$ donnant un contr\^ole de
$u_\ep$ en $\ep^4$. 
La fonction $\phi_\ep =
\max (\phi , \Re (h-u_\ep)+\ep^3)$ est pluriharmonique sur un petit
voisinage de $S$, mais est \'egale \`a $\phi $ en-dehors de $S_\ep$ du fait
de la croissance de $\phi - \Re h$.

\subsection{R\'esum\'e de la d\'emonstration du cas immerg\'e}

Dans le cas o\`u $S$ est seulement immerg\'ee,
la pr\'esence des points doubles rend impossible, sans travail
suppl\'ementaire,  le prolongement $\db$-plat. Par ailleurs, 
m\^eme si ce probl\`eme est r\'esolu,
la diff\'erence $\phi- \Re h $ ne cro\^{\i}t plus
automatiquement au voisinage des  points doubles.\\

Les deux \'etapes de la d\'emonstration du cas plong\'e n\'ecessitent
chacune une r\'esolution du $\db$: d'abord pour cr\'eer une forme
ferm\'ee de type $(1,1)$, positive \`a support hors de $S$, ensuite
pour approcher son
support par une hypersurface. Notre premier lemme
 condense les deux \'etapes en une
seule: la platitude du $\db$ de l'extension $h$ suffit
\`a repousser hors de $S$ une hypersurface passant par un point
du compl\'ementaire.

On conserve n\'eanmoins une partie de la deuxi\`eme \'etape pr\'ec\'edente:
il nous faut
modifier $\om$ pour l'annuler au voisinage des points doubles
tout en conservant le caract\`ere lagrangien de  $S$. 
Si $\phi$ est un potentiel de la forme modifi\'ee,
on construit
une fonction $h$ globale $\db$-plate sur $S$, dont
la partie r\'eelle est une extension de $\phi_{|S}$
hors de $S$,
 avec $h$ holomorphe pr\`es des points doubles.

La croissance en $d^2(.,S)$ de  $\phi -\Re h$ est encore
assur\'ee l\`a  o\`u $\phi$ est strictement plurisousharmonique, mais
plus au voisinage des points doubles. Nous modifions
alors $\phi$ sans changer $\phi_{|S}$ pour la rendre strictement 
plurisousharmonique \`a l'endroit crucial: la trace sur $S$ de la fronti\`ere
du support de $\omega$.
Cette ultime modification est technique,
et se fonde essentiellement sur la convexit\'e
polynomiale de la r\'eunion du cylindre unit\'e $\{|\Re z| \le 1\}$
 et de l'espace totalement r\'eel $\Rr^n$, ou de la boule et de $\Rr^n$,
 dont l'\'etude remonte au moins \`a Smirnov et Chirka ([Ch,Sm], voir aussi
[Bo]). Plus pr\'ecis\'ement,
 nous explicitons une fonction
 plurisousharmonique positive s'annulant exactement sur une de ces deux figures.\\

\textbf{Remerciements.} Je tiens \`a remercier Julien Duval, qui m'a donn\'e
ce sujet de recherche,  pour son soutien
continu durant ce travail, pour la richesse
de ses intuitions g\'eom\'etriques, et enfin pour 
ses patientes relectures de cet article.


\section{Preuve
du th\'eor\`eme}
\subsection{Un lemme pour la
convexit\'e rationnelle}

Dans toute la suite, $S$ d\'esignera une sous-vari\'et\'e $C^\infty$ 
immerg\'ee \`a points doubles ordinaires, totalement
r\'eelle de dimension moiti\'e, compacte et sans bord.
L'\'enonc\'e suivant permet, dans le cas plong\'e,
 de construire directement une hypersurface
passant par un point du compl\'ementaire de $S$, \`a partir
d'une fonction $\db$-plate sur $S$:\\

\textit{
Soit $\phi$  strictement plurisousharmonique
 $C^\infty$ sur $\Cc^n$ et $S$ plong\'ee.
S'il existe $h$
une fonction $C^\infty$ v\'erifiant les propri\'et\'es suivantes :\\
                 \hspace*{.5cm}
$\bullet \ \db h = O(d^{m}(., S))$ o\`u $m=\frac{1}{2}(3n+5)$,\\
                \hspace*{.5cm}
$\bullet \ |h|=e^\phi $ \`a l'ordre 1 sur S,\\
alors S est rationnellement convexe.
        }\\

         L'existence de $h$ est naturelle, car elle implique que  $S$ est
        lagrangienne pour $dd^c\phi$.

    Cet \'enonc\'e n'est pas suffisant dans le cas immerg\'e,
        et sera la cons\'equence du lemme suivant:

                                \begin{lemma}
        Soit $\phi$ plurisousharmonique  $C^\infty$ 
        sur $\Cc^n$, et $h$ une fonction
        $C^\infty$  sur $\Cc^n$ telle que  $ |h|\le e^\phi$, avec :\\ 
        \hspace*{.5cm}
        $\bullet \ \db h = O(d^{\frac{3n+5}{2}}(.,S))$\\
        \hspace*{.5cm}
        $\bullet \ |h|= e^\phi$ \`a l'ordre 1 sur S, 
        et $X= \{|h|= e^\phi\}$ compact,\\
        \hspace*{.5cm}
        $\bullet$ pour tout x de X,
         une des deux  conditions suivantes est remplie :\\
                 \hspace*{1cm}
        1) $h$ est holomorphe sur un voisinage de $x$ \\
                \hspace*{1cm}
        2) x est un point r\'egulier de $S$,
         et $\phi$ est strictement plurisousharmonique
      en x.\\
        Alors X est rationnellement convexe.
                                \end{lemma}

Ce lemme nous permet de travailler avec une fonction $\phi$ seulement
 plurisousharmonique, 
 avec les  conditions, entre autres,  qu'elle soit strictement 
plurisousharmonique aux points r\'eguliers de $S$ pr\`es desquels
$h$ n'est pas holomorphe, et que $h$ soit 
holomorphe pr\`es des points doubles.

Le premier \'enonc\'e est un corollaire de ce lemme, car 
\textit{dans le cas o\`u $\phi$ est strictement plurisousharmonique, 
       $ \phi - \log |h|$ cro\^\i t automatiquement 
        en $d^2(.,S)$}, si bien qu'on peut prendre $X=S$.  
        En effet,
        puisque $\db h = O(d^2(., S))$, $dd^c h$ est nulle sur $S$.
        Donc $\psi= \phi - \log|h| $
        est strictement plurisousharmonique  sur un voisinage de $S$,
        et s'annule \`a l'ordre 1 sur $S$.
        Gr\^ace \`a  un param\'etrage local de $S$ par
        $\Rr^n$ d'extension $\db$-plate,
         donn\'e
        par le lemme 2, on peut se ramener au cas o\`u $S=\Rr^n$.
        Si $x=\Re z$ et $y=\Im z$, on a donc
        $$
         \psi (x+iy) = \sum_{j,k} \frac{\partial^2 \psi} {\partial y_j \partial
         y_k}(x) y_j y_k + O(|y|^3).
        $$
        Par ailleurs, la hessienne r\'eelle sur
        $\Rr^n$ de $\psi$ co\"{\i}ncide  \`a un facteur positif pr\`es avec  sa
        hessienne complexe qui est d\'efinie positive.
        Donc
        $\psi(z) \ge  A |y|^2$ pour un $A>0$.\\

\begin{preuve}[du lemme 1] Soit   $B$ une boule contenant $X$, 
                et $x_0 \in B\setminus X$.  Il nous faut trouver une
                hypersurface \'evitant $X$ et passant par $x_0$. 
                 Quitte \`a multiplier 
                $h$ par une fonction valant $1$ sur un voisinage de $X$ mais
                nulle pr\`es de $x_0$, on peut supposer que $h$ s'annule sur
                une petite boule centr\'ee en $x_0$.
                Soit $u = h^N$, o\`u
                $N$ est un entier naturel qu'on fera tendre vers l'infini.
                D\'efinissons le poids plurisousharmonique $\psi$ 
                par $\psi(z) = 2N\phi(z) +
                2n\ln |z-x_0| + |z|^2$.
                Gr\^ace au th\'eor\`eme de H\"ormander (cf. [H\"o]), nous
                r\'esolvons l'\'equation $\db v = \db u$ sur la boule,
                avec l'estimation $\|v \|^2_\psi \le C\| \db u\|^2_\psi$,
                et $C$ ne d\'ependant que de la boule. 
                L'in\'egalit\'e s'\'ecrit :
                $$
                \int_B \frac {|v|^2 e^{-2N\phi -|z|^2}}
                {|z-x_0|^{2n}} d\lambda \le C
                \int_B \frac {|\db u|^2 e^{-2N\phi -|z|^2}} {|z-x_0|^{2n}}
                d\lambda. $$
                La premi\`ere int\'egrale converge, donc $v(x_0) = 0$, et
                l'hypersurface d\'efinie par $\{ f= v-u= 0\}$ passe
                par $x_0$.  Pour v\'erifier qu'elle ne rencontre pas $X$,
                il suffit de
                trouver $N$ assez grand, tel que 
                $|v(z)| \le \frac{1}{2}e^{N\phi(z)}$
                pour tout $z \in X$.
                Par le lemme de H\"{o}rmander-Wermer ([H\"{o},We]), on a, pour
                $z \in X$ et $\ep >0$,
                                \begin {eqnarray*}
                |v(z)|^2 & \le &
                        C_1\left( \ep^2 \|\db v\|^2_{L^{\infty}B(z,\ep)} +
                        \ep^{-2n} \|v\|^2_{L^2(B(z,\ep))}\right)\\
                 & \le & C_1\left( N
                ^2\ep^2 \|\db h\
                h^{N-1}\|^2_{L^\infty(B(z,\ep))} + \ep^{-2n} \|v e^{-
                \frac{\psi}{2}} 
              e^{ \frac{\psi}{2} } \|^2_{L^2(B(z,\ep))} \right).\\
                                \end{eqnarray*}

                Sur la boule $ B(z,\ep)$, on a $|h(x)| 
                \le e^{\phi(z)+c\ep}$, o\`u $c=\|\phi\|_{C^1(B)}$.  D'o\`u

                                \begin{eqnarray*}
        \|\db h\ h^{N-1}\|^2_{L^\infty(B(z,\ep))} & \le & C_2
        e^{2N\phi (z)}e^{2Nc\ep} \|\db
        h\|^2_{L^\infty(B(z,\ep))},\\
        \|v e^{- \frac{\psi}{2}} e^{\frac{\psi}{2} }
        \|^2_{L^2(B(z,\ep))} & \le & C_3 e^{2N\phi
        (z)}e^{2Nc\ep} \|v e^{- \frac{\psi}{2}} \|^2_{L^2(B(z,\ep))}.\\
                                 \end{eqnarray*}


Choisissons $\ep = \frac{1}{N}$.  L'in\'egalit\'e devient
$$|v(z)|^2 \le C_4 e^{2N\phi (z)}
\left ( \|\db h\|^2_{L^\infty(B(z,\frac{1}{N}))} + N^{2n} \|v e^{-
\frac{\psi}{2}} \|^2_{L^2(B(z,\frac{1}{N}))}\right). $$

On a $\db h(x) = O( d^{\frac{3n+5}{2}}(x, S))$, tandis que  $\db h(x)\equiv 0$ 
au voisinage de $z \in X$ si $z$ n'est pas un point r\'egulier de $S$.  Dans 
tous les cas,
$\|\db h\|^2_{L^\infty(B(z,\frac{1}{N}))} \le C_5
\frac{1}{N^{3n+5}} \underset{N\to \infty}{\to} 0 $. \\
Pour le second terme de la parenth\`ese, on a: 
                \begin{eqnarray*}
 N^{2n} \|v e^{- \frac{\psi}{2}} \|^2 _{L^2(   B(z,\frac{1}{N})   ) }
 & \le & N^{2n} \|v\|_\psi^2 \\
 & \le & C \ N^{2n} \|\db u\|^2_\psi\\
 & \le & C' N^{2n+2} \int_{B} |\db h|^2\left(|h| e^{-\phi}\right)^{2N}.\\
                \end{eqnarray*}
Soit $\alpha >0$.
Nous allons couper l'int\'egrale en deux, l'une ayant pour domaine le
voisinage tubulaire de $S$ de taille $1/{N^\alpha}$, l'autre le compl\'ementaire
dans la boule :
                \begin{eqnarray*}
N^{2n+2}\int_{Tub(S,\frac{1}{N^\alpha})} |\db h|^2 \left(|h|
e^{-\phi}\right)^{2N} & \le & \ N^{2n+2}\int_{Tub(S,\frac{1}{N^\alpha})}
|\db h|^2 \\
                      & \le & c N^{2n+2}\left (\frac{1}{N^{\alpha}}
\right)^{3n+5} \left
(\frac{1}{N^{\alpha}}\right)^{n}. \\
                \end{eqnarray*}


Il suffit donc, pour que cette int\'egrale tende 
vers $0$ quand $N$ tend vers l'infini,
de prendre $$ \alpha >\frac {2n+2}{4n+5}.$$

Pour majorer  la deuxi\`eme partie de
l'int\'egrale, rappelons qu'au voisinage d'un  point r\'egulier de $S$ o\`u $\phi$ est
strictement plurisousharmonique, on a l'estimation 
$\phi -\log |h| \ge ad^2(z,S)$, o\`u $a$ est
une constante strictement positive. 
Par ailleurs, d'apr\`es l'alternative donn\'ee par 
les hypoth\`eses du lemme, un point 
$x \in B$ ou bien se trouve dans une boule centr\'ee en un point de $X$
et sur laquelle $h$ est holomorphe, ou bien appartient \`a une r\'egion
proche de $S$ o\`u la minoration pr\'ec\'edente est valable, ou bien 
enfin se situe hors   de ces deux ensembles, et dans ce dernier cas 
$|h|e^{-\phi} \le k<1$ pr\`es de $x$ ($k$ est une constante uniforme). Par cons\'equent
\begin{eqnarray*}
         N^{2n+2}\int_{B\setminus Tub(S,\frac{1}{N^{\alpha}})}
             |\db h|^2 \left(|h| e^{-\phi}\right)^{2N}
  & \le & c N^{2n+2}\int_{B\setminus Tub(S,\frac{1}{N^{\alpha}})}
             \left ( 1 - a d^2 (x,S)\right ) ^{2N} \\
           &  + &c' N^{2n+2} k^{2N}.\\
\end{eqnarray*}
La derni\`ere int\'egrale est inf\'erieure \`a
 $$c' N^{2n+2} \left ( 1 - \frac{a}{N ^{2\alpha}}\right ) ^{2N}
 \sim c' N^{2n+2}\exp( -2a N ^{1-2\alpha}),
$$
tandis que $N^{2n+2} k^{2N}\to 0$.
 Il suffit donc de
choisir, pour que cette partie de l'int\'egrale 
converge vers $0$ quand $N$ tend vers
l'infini :  $$\alpha <1/2.$$ Cette condition est compatible avec
la pr\'ec\'edente.  On peut alors prendre $N$
assez grand ind\'ependamment de $z\in S$, de sorte que $|v(z)| \le
\frac{1}{2}e^{N\phi(z)}$.  On aura $|f(z)| =|e^{N\phi(z)} - v(z)| \ge
\frac{1}{2}e^{N\phi(z)} \not= 0$, et l'hypersurface \'evite $X$.
\end{preuve}


  La d\'emonstration du th\'eor\`eme se ram\`ene maintenant \`a
  satisfaire les conditions
  du lemme 1. Plus pr\'ecis\'ement, nous allons 
  construire des ensembles $X$ contenus dans 
  la r\'eunion de $S$ avec des voisinages  
  des points doubles, asymptotiquement proches de boules aussi petites
  que l'on veut. $S$ sera alors rationnellement convexe.
  
  Pour le confort du lecteur, nous supposerons dor\'enavant que
  $S$  poss\`ede un unique point double en l'origine.
  Dans le paragraphe 2.2, nous annulons la forme de K\"{a}hler
  sur une petite boule centr\'ee sur $0$, 
  tout en laissant $S$ lagrangienne.
  Ensuite, dans le paragraphe 2.3, 
   nous exhibons le prolongement  $h$
  $\db$-plat sur $S$ et holomorphe au voisinage de l'origine.
  La construction de $h$ est globale, mais la d\'emonstration
  est identique s'il y a plus d'un point
  double.
  Enfin, dans le paragraphe 2.4,
  nous relevons  $\phi$ en une fonction 
  strictement plurisousharmonique sans changer $\phi_{|S}$, pr\`es de $0$
 et pr\'ecis\'ement \`a l'endroit o\`u $\phi$ n'est  pas encore
  strictement plurisousharmonique mais o\`u $h$ perd son caract\`ere holomorphe.

\subsection{ Annulation  de $\om$ au voisinage du point double}
        Quitte \`a retrancher la partie r\'eelle d'un polyn\^ome 
        quadratique \`a un potentiel $\phi$
        de $\omega$,
        on peut supposer que $\phi$ est strictement
        convexe au voisinage de l'origine, et que  $0$ est
        un minimum local strict pour $\phi$.
        Par un changement lin\'eaire
        de coordonn\'ees, on 
        se ram\`ene au cas o\`u 
        $\phi(z) = |z|^2 + O(|z|^3)$.
         La proposition suivante
        annule la forme $\omega$ sur une petite 
        boule contenant $0$,
        mais sans  changer $\om$ hors d'une boule plus grande, tout
        en conservant le caract\`ere lagrangien de $S$. De plus
        le domaine d'incertitude concernant
        la stricte positivit\'e de la forme transform\'ee
        doit \^etre aussi petit qu'on veut.
         La taille de ce no man's land
        sera en effet prescrite par l'angle form\'e
        par les deux branches de $S$ en l'origine.

                        \begin{prop}
        Soit S et  $\om$ comme ci-dessus.
         Alors pour tout couple de boules assez petites et centr\'ees sur le point
         double, il
        existe une $(1,1)$-forme $\tilde{\om}$ positive telle que:\\
                \hspace*{.5cm}
           $\bullet $ S est lagrangienne pour
        $\tilde{\om}$,\\
                \hspace*{.5cm}
         $\bullet \ \tilde{\om}$ est nulle sur
        la plus petite boule,\\
                \hspace*{.5cm}
        $\bullet \ \tilde{\om} = \om$ hors de la seconde boule.
                                \end{prop}

\begin{preuve}
Nous allons construire, pour tout $\ep >0$ assez petit,
 une fonction $\zeta_{\ep}$ telle que
$j^*d^c \zeta_{\ep} = j^*d^c \phi$ sur une
boule $b_\ep$ de taille $\ep$, 
de support contenu dans  une boule un 
peu plus grande $b^+_\ep$, avec $j^*dd^c \zeta_{\ep} =0 $, et telle que les d\'eriv\'ees d'ordre
$k=0,1,2$ soient contr\^ol\'ees par $|z|^{3-k}$, uniform\'ement
en $\ep$.  Dans ce cas la fonction 
$\phi_\ep = \phi - \zeta_\ep$ 
reste strictement plurisousharmonique,
v\'erifie $j^*d^c \phi_\ep=0$ sur $b_\ep$, $j^*dd^c \phi_{\ep} =0 $ et 
$\phi_\ep(z) = |z|^2 + O(|z|^3)$.
        Soit alors $\rho$ une fonction positive $C^\infty$, convexe et
        croissante d\'efinie sur $\Rr$, nulle pour $t$ assez petit et
        valant $t-c$ pour $t$ un peu plus grand.
        Si les param\`etres de $\rho$ sont choisis de fa\c{c}on ad\'equate
          en fonction des deux petites boules,
        la  (1,1)-forme  positive 
        $\tilde{\om} = dd^c (\rho \circ \phi_\ep)$ est 
  nulle sur une boule $b^-_\ep$ proche de $b_\ep$, est
  \'egale \`a $dd^c \phi_\ep$ hors de $b_\ep$ et donc  \`a $dd^c \phi$ 
  hors de $b^+_\ep$.  Sachant que 
 $j^*d^c (\rho \circ \phi_\ep) =0$ sur $b_\ep$,
 on v\'erifie alors ais\'ement 
 que $S$ reste lagrangienne pour $\tilde{\omega}$.\\

Nous d\'emontrons maintenant l'existence de $\zeta_\ep$.
Nous construisons en fait  $\zeta_1$ et  $\zeta_2$,
deux fonctions r\'ealisant les propri\'et\'es requises par $\zeta_\ep$, mais
pour $S=S_1$ et $S=S_2$, o\`u $S_1$ et $S_2$ sont les deux branches de $S$ 
en $0$. Si  $\gamma_1$ est
une fonction $C^\infty$ d\'efinie sur la sph\`ere unit\'e
$S^{2n-1}$, \'egale \`a 1 sur un
voisinage de l'intersection de $T_0 S_1$ avec $S^{2n-1}$,
 et nulle sur un voisinage de
l'intersection de $T_0 S_2$ avec la sph\`ere, alors les d\'eriv\'ees d'ordre
$k= 0,1$ ou $2$ de la  fonction
$\tilde\zeta_1 (z)= \gamma_1 (\frac {z}{|z|})\zeta_1 (z)$
 sont encore en $O(|z|^{3-k})$
uniform\'ement en $\ep$. La fonction 
$\zeta_\ep = \tilde\zeta_1 + \tilde\zeta_2$
est celle que nous recherchions ($\tilde\zeta_{2}$ est construite
comme $\tilde\zeta_{1}$, mutatis mutandis). Dans le reste de ce paragraphe, nous pouvons
donc supposer que $S$ est plong\'ee en 0.

 La forme $j^*d^c \phi$ est ferm\'ee,
  donc poss\`ede une primitive $\lam$
  d\'efinie pr\`es
  de 0  sur $S$. Par ailleurs on peut 
  la choisir telle que $\lam(z) = O(|z|^3)$, 
  car $j^*d^c \phi (z) = O(|z|^2)$.
  En effet, dans nos coordonn\'ees, $P = T_0 S $ est lagrangien
  pour la forme standard $\om_0 = dd^c |z|^2$. Il est facile
  de v\'erifier que $j_P ^* d^c |z|^2=0$, o\`u $j_P $ est l'injection
  de $P$ dans $\Cc^n$. L'estimation s'en d\'eduit.

  On coupe cette fonction en $\lam_\ep (z) = \chi (\frac{z}{\ep}) \lam(z)$
  par une fonction plateau
  valant 1 sur la boule de rayon $\ep$, et de support dans une boule
  un peu plus grande.
   On v\'erifie ais\'ement que
  les d\'eriv\'ees d'ordre $k =0,1,2,3$ 
  sont contr\^ol\'ees par $|z|^{3-k}$, et cela
  uniform\'ement en $\ep$.  Nous voulons maintenant trouver
  une fonction $\zeta_\ep$ telle que $j^*d^c \zeta_\ep = d\lam_\ep$.
  Il suffit de prendre $\zeta_\ep $ nulle 
  sur $S$ et de sp\'ecifier ses d\'eriv\'ees
  le long des directions transverses \`a $TS$, 
  ce qui est possible car $S$ est
  totalement r\'eelle: si $S = \Rr^n$ et
  $d\lam_\ep = \sum_i \alpha_i (x) dx_i$, la fonction
  $\zeta_\ep(z) = \sum_i \alpha_i (x) \cdot y_i $
  convient, et ses d\'eriv\'ees d'ordre $k=0,1,2$ sont estim\'ees par
  $|z|^{3-k}$ uniform\'ement  en $\ep$.
 Dans le cas g\'en\'eral, on se ram\`ene au cas pr\'ec\'edent par
 le lemme classique ci-dessous. \end{preuve}


                        \begin{lemma}
        Soit $S$ une sous-vari\'et\'e lisse de
        $\Cc^{n}$ passant par $0$, avec $T_0 S $   totalement r\'eel.
        Alors pour tout $m\in \Nn$, il existe un
        diff\'eomorphisme local $\psi$ fixant 0
        tel que $\psi (S) \subset \Rr^{n}$ et $\db
        \psi (z) = O(d^m(z,S)) $.
                        \end{lemma}
\begin{preuve}
 On peut supposer que $T_0 S = \Rr^n$. 
 Il existe $\phi: S \to \Rr^n$  un diff\'eomorphisme 
 local pr\`es de $0$ tangent \`a l'identit\'e. 
 Posons $\phi = (\phi_1, \cdots ,\phi_n)$.
 D'apr\`es [H\"{o}, We], on peut prolonger chaque $\phi_k$ en $\psi_k$ en
 dehors de $S$, avec
 $\psi_{k|S}= \phi_{k|S}$, et
 $\db \psi_k = O(d^m(.,S))$.
 Donc $D_0 \psi$ est $\Cc$-lin\'eaire et 
 l'identit\'e sur $\Rr^n$. On en d\'eduit que $\psi$
 est tangent \`a l'identit\'e et donc qu'elle est  un diff\'eomorphisme local. 
 \end{preuve}
 
\subsection{Construction d'une extension
            $\bar\partial$-plate}
           Soit $\phi$ un potentiel de $\tilde \omega$, la forme fournie par
           la proposition 1 pour un couple de petites boules centr\'ees sur
           le point double. Rappelons que $\phi$ est nulle sur la petite boule
           et strictement plurisousharmonique hors de la seconde. Il s'agit 
           maintenant de construire une fonction $h$ $\db$-plate sur $S$
           et telle que $\phi - \log |h|$ s'annule \`a l'ordre 1 sur $S$.
           Ces deux contraintes entra\^{i}nent $j^*( d^c \phi - d (\arg h)) =0$. 
           Autrement dit, la 1-forme ferm\'ee  $j^* d^c \phi$ doit avoir
           des p\'eriodes enti\`eres (ou plus g\'en\'eralement rationnelles, 
           quitte \`a multiplier $\phi$ par un entier). 
           Comme dans  [Du,Si], il nous faut perturber 
           $\phi$ pour r\'ealiser cette condition.

        Soit $\tilde {S}$ un ouvert lisse de $\Cc^n$  tel que
        $S$ soit un retract par d\'eformation de $\tilde {S}$.
        On a alors $ H_1(\tilde {S}, \Zz)\simeq H_1(S, \Zz)$.
        Soit
        $\gamma_1, \cdots,\gamma_p$ une base de $ H_1(\tilde {S}, \Zz)$,
        qu'on peut supposer support\'ee par $S$.
        D'apr\`es le th\'eor\`eme de De Rham, il existe des $1-$formes ferm\'ees
        $\beta_1, \cdots,\beta_p$  sur $\tilde {S}$,
        telles que $\int_{\gamma_i} \beta _j = \delta _{ij}$.  On peut
        choisir ces formes nulles sur une boule $b$ contenant
        le point double de $S$, et telle que $2b \Subset \tilde{S}$.
        Il suffit pour cela de prendre 
        une primitive $f_j$ de chaque $\beta_j$ sur $2b$, et $\chi$ une fonction
        plateau nulle sur $b$ et valant $1$
        hors de $2b$. La 1-forme $\chi \beta_j + f_j d\chi$ est
        dans la m\^eme classe que $\beta_j$ et s'annule sur $b$.
        Il est alors possible de trouver
        $\psi_1, \cdots, \psi_p$ \`a support compact dans
        $\Cc^n$, v\'erifiant $j^*d^c \psi_j = j^*\beta_j$,
        et $\psi_j \equiv 0$ sur $b$.
        En effet, on peut fixer $\psi_j \equiv 0$ sur
        $S \cup  b$, et sp\'ecifier ses d\'eriv\'ees
        dans les directions contenues dans $ iT_zS$
         (cf. d\'emonstration de la proposition 1).
        
        Posons maintenant
        $\phi_\lambda = \phi + \lambda_1\psi_1 + \cdots + \lambda_p \psi_p$.
        On a $\phi_\lambda = \phi$ sur $S \cup  b$.
        Pour $\lambda = (\lambda_1 , \cdots , \lambda_p)$ assez petit,
        $\phi_\lambda$ est encore strictement plurisousharmonique en dehors de
        la plus grande des deux petites boules et plurisousharmonique sur $b$.
        On peut de plus trouver $M$ entier et choisir $\lambda$, tels que
        $\int_{\gamma_i}j^*d^c \phi_\lambda \in 2 \pi \Zz/M$ pour $1\le i\le p$.
        Posons $\tilde\phi= M\phi_\lambda$. La forme
        $j^*d^c \tilde\phi$ est ferm\'ee sur $S$ et ses p\'eriodes sont des
        multiples de $2\pi$.  Il existe donc une fonction $C^\infty$
        $\mu :  S \rightarrow \Rr/{2\pi \Zz}$, 
                telle que $j^*d^c \tilde\phi = d\mu$. On peut \'evidemment choisir
                $\mu$ nulle sur l'intersection de $S$ avec la petite boule. 
        D'apr\`es [H\"{o},We], pour tout entier $m$,
         il existe une fonction $h$ d\'efinie sur $\Cc^n$, telle
        que $h_{|S} = \exp(\tilde\phi + i\mu)_{|S}$ avec
        $\db h (z) = O(d^m (z,S))$. Automatiquement,
        $ \tilde\phi - \log|h|$ s'annule \`a l'ordre 1  sur $S$.
        Pour ne pas alourdir les notations, nous continuons \`a nommer $\phi$
        la nouvelle fonction $\tilde \phi$.
        
        Dans le paragraphe suivant, nous rendons $\phi$ strictement
        plurisousharmonique sur un petit ouvert contenant
        l'intersection de $S$ avec le no man's land o\`u l'on ne sait rien
        sur la stricte positivit\'e de $dd^c \phi$, et ce sans
        changer $\phi_{|S}$. 
        On peut alors prolonger $h$ par $1$ sur un voisinage du point double, 
        de sorte que 
         $X= \{e^\phi= |h|\}$ soit un compact v\'erifiant les conditions
        du lemme 1. $\Box$

\subsection{Rel\`evement de $\phi$}

        Comme il a \'et\'e annonc\'e pr\'ec\'edemment, 
        il est n\'ecessaire  maintenant de relever
         $\phi$  pour la rendre strictement plurisousharmonique 
         sur un petit voisinage de $\partial (\text{support}\  dd^c \phi)
         \cap S$, i.e l\`a o\`u $h_{|S}$ ne peut plus \^etre prolong\'ee de
         fa\c{c}on holomorphe et o\`u $\phi$ n'est pas encore strictement
         plurisousharmonique. On peut supposer que  $T_0 S_1 = \Rr^n$.
        
           Nous quantifions la latitude donn\'ee par la proposition 1
        sur le choix des petites boules: prenons $\delta$, telle que  
        $T_0S_2 \cap B(0,1+4\delta) \subset \{|x| \le 1\}$ (cf. Fig. 1; rappelons
        que $x=\Re z$ et $y=\Im z$). 
        On peut alors choisir les deux boules de
        sorte qu'il soit suffisant de relever $\phi$ entre
        $B(0, (1+2\delta)\ep )$ et $B(0, (1+3\delta)\ep)$ au voisinage
        de $S_1$ et de $S_2$.
         
        Nous construisons d'abord une fonction $\chi_\ep$ plurisousharmonique
        sur un voisinage
        tubulaire de $S_1$ de taille $\sim \ep$,
        nulle sur $S_1$ et sur le cylindre $\{|x|\le (1+\delta)\ep\}$, et
        croissant entre les deux  boules en $d^2(z,S_1)$.
        Nous avons ensuite besoin de la prolonger hors du tube par
        une fonction $\psi$ plurisousharmonique nulle sur $S_2$, au moins
        sur la boule hors de laquelle $\phi$ est strictement
        plurisousharmonique. 
        A l'ext\'erieur de cette  boule, on pourra couper $\psi$ en profitant
        de la stricte plurisousharmonicit\'e de $\phi$.
         Le lemme suivant fournit cette fonction, et
        prouve la convexit\'e polynomiale de la r\'eunion du
        cylindre et de $\Rr^n$:

                        \begin{lemma}
        Il existe une fonction $\psi$ $C^\infty$, 
        plurisousharmonique et positive sur $\Cc^n$,
        nulle exactement sur la r\'eunion du  
        cylindre unit\'e $\{|x| \le 1\}$ avec $\Rr^n$.
                        \end{lemma}
\begin{preuve}
        Soit $\alpha$ une fonction $C^\infty$ de $\Rr^+$
        dans $\Rr^+$, convexe, croissante, nulle sur
        $[0,1]$ et strictement positive sur $]1, \infty]$. 
        La fonction  
        $\psi_1 (z) = \alpha (|x|) y^2 $
        est  sous-harmonique et cl\^ot la d\'emonstration pour $n=1$.

        Pour $n\ge 2$, posons
        $$ \psi (z) = \int_{S^{n-1}} \psi_1( <z,X>) d\sigma (X),$$
        o\`u $d\sigma $ est la mesure normalis\'ee invariante par  $O_n(\Rr)$
           sur la sph\`ere unit\'e
        $S^{n-1}$ de $\Rr^{n}$,
         et $<.,.>$ le produit scalaire hermitien
        canonique de $\Cc^n$.
        Comme int\'egrale de fonctions plurisousharmoniques,
         $\psi $ est plurisousharmonique.
        Pour $z$ appartenant \`a $\Rr^n$ ou au  cylindre $\{|x|\le 1\}$, $\psi(z) = 0$.
        V\'erifions maintenant que $\psi$ est strictement
        positive en dehors de cette figure. Pour cela, remarquons
        que l'application d\'efinie sur $S^{n-1}$ par 
        $X \mapsto  \alpha (|<x,X>|) <y,X>^2 $ pour
        $|x| >1$ et $y\not= 0$ s'annule sur 
        une calotte ne couvrant pas $S^{n-1}$ tout
        enti\`ere et sur l'intersection d'un hyperplan
        avec la sph\`ere. Le compl\'ementaire de la r\'eunion de ces deux ensembles est
        un ouvert non vide de $S^{n-1}$, et donc $\psi(z) >0$.
                        \end{preuve}

                 \textit{Remarque}: On peut
         prouver de la m\^eme fa\c{c}on la convexit\'e
         polynomiale du compact $L$ form\'e par la r\'eunion de la boule
         unit\'e de $\Cc^n$ avec la boule de rayon 2 dans $\Rr^n$.
         En effet, soit dans $\Cc$ la fonction de
         Green $G_1$ de $\{|z|<1\} \cup [-2,2]$,
         avec p\^ole logarithmique \`a l'infini. 
                Comme le compact est r\'egulier pour le probl\`eme
                de Dirichlet, la fonction 
                est sousharmonique continue  sur $\Cc$, 
                positive et nulle 
                pr\'ecis\'ement sur le compact.
                La fonction  $G$ d\'efinie par une int\'egrale analogue \`a
         celle de $\psi$ dans le lemme ci-dessus
         est plurisousharmonique et positive, nulle pr\'ecis\'ement 
                sur le compact $L$.  En particulier
                ses sous-niveaux  produisent directement
         une base de voisinages pseudoconvexes de la figure $L$ (comparer
         \`a [Bo]).

\unitlength=1cm
\begin{picture}(14,10)
 
 \put(9,1.1){\vector(0,1){.8}  }        
 \put(9.2,1.5){$\sqrt{\frac{c}{2}}\ep$}        
 \put(3,2){\dashbox{.1}(6,0)  }         
 
 \put(10,.5){$\Rr^n=T_0S_1$}
 \put(1.8,7.7){$T_0S_2$}
 \put(8.2,1.2){$S_1$}
 \put(0,9){$i\Rr^n$}

 \put(3,1){\dashbox{.1}(0,7)  }
 \put(2.8,.5){\small{$\epsilon$}}
 
 \put(5,1){\dashbox{.1}(0,1)  }
 \put(4.7,.5){\small{$(1+\delta)\epsilon$}}
 
 \qbezier(3.3,5)(3,1.4)(6,1.5)                  
 \put(3.3,5){\line (0,1){3}}
 \put(6,1.5){\line (1,0){2.2}}
 \put(3.7, 7.4){\vector(-1,0){.4}}
 \put(4, 7.4){$\{ \psi_\ep  = 0 \} $}
 
 \qbezier[200](7,1)(6.9,6.9)(1,7)               
 \put(6.7,.5){\small{$(1+4\delta)\epsilon$}}

  \thicklines
  \put(1,1){\vector(1,0){9.2}  }         
  \put(1,1){\vector(0,1){8}  }
  
  \put(1,1){\line  (1,4){1.8}}          
  \qbezier(1,1)(6,1.1)(8,1.3)           
  
\end{picture}
        \begin{center}
        \small{\textbf{Fig. 1 } }
        \end{center}
        
\vspace{1cm}

          Nous construisons maintenant la fonction $\chi_\ep$. Soit pour cela 
           $\chi$ une fonction $C^\infty$ n\'egative, telle que
                $ \chi(t) = -1 $ pour $t \le  1+\delta$, et
                $ \chi(t) = 0 $ pour $t \ge 1+2\delta$. Posons
        $$ \chi_\ep (z) = \max \left( d^2_1(z) +
                c \ep^2
                \chi(\frac{|x|}{\ep}), 0 \right) ,$$
        o\`u $d_1 $ est la distance \`a $S_1$.
        Puisque
        $ \|\ep^2 \chi(\frac{|x|}{\ep}) \|_ {C^2} =
        O(1) $, et que $d^2_1$ est
        strictement plurisousharmonique sur un voisinage fixe
        de $S_1$, $ d^2_1(z) + c \ep^2
        \chi(\frac{|x|}{\ep})$ l'est encore pour $\ep$ et $c$ assez
        petits ($c \sim \delta^{2}$), et $\chi_\ep $ est plurisousharmonique. 
        Notons que si  $\ep$ est assez petit, pour 
        $|y| \le \sqrt{\frac{c}{2}}\ep$ et $|x|\le (1+\delta) \ep $,
        $\chi_\ep \equiv 0$ (car $d_1(z)= |y| + O(|y|^2)$, 
        tandis que  pour  $  (1+2\delta) \ep \le |x| $,
        $ \chi_ \ep \equiv d_1^2$.

        Nous prolongeons maintenant cette fonction
         de fa\c{c}on plurisousharmonique
        en dehors du tube $|y| = \sqrt{\frac{c}{2}}\ep$.
         Nous utilisons pour cela 
        la fonction  $\psi$ donn\'ee par le lemme 3.
        On d\'emontre facilement qu'il existe une constante $k>1$, 
        telle que pour
        $(1+\delta) \le |x| \le (1+4\delta)$, on ait l'estimation 
        $k^{-1}|y|^2 < \psi(z) < k|y|^2 $. La fonction
        $$\psi_\ep(z) = \psi (\frac{z}{\ep})- \frac{c}{4k}$$
        v\'erifie, pour $(1+\delta)\ep \le |x| \le (1+4\delta)\ep$, et $\ep $ 
        assez petit, $\psi_\ep > \frac{c}{4k} $ sur 
        $|y|=  \sqrt{\frac{c}{2}}\ep$ et 
        $\psi_\ep < 0 $ sur $S_1$ (cf. Fig. 1).
        Or $\chi_\ep \le
        \frac{c}{2}\ep^2$ sur le bord du tube. 
        Donc la fonction $\max (\psi_\ep, \chi_\ep)$ 
        est plurisousharmonique, nulle sur $S_1$ et sur
         $S_2\cap B(0, (1+4\delta)\ep)$, 
        et est strictement plurisousharmonique sur un voisinage de $S_1$
        entre  $B(0, (1+\delta)\ep)$ et $B(0, (1+4\delta)\ep)$.
         Si $\gamma$ est une
        fonction plateau valant $1$
        sur $B(0,(1+4\delta)\ep)$ et
        nulle peu apr\`es, 
        la fonction $\phi + a\gamma \max (\psi_\ep, \chi_\ep)$
        est plurisousharmonique sur $\Cc^n$ pour $a$ assez petit
         et reste \'egale \`a $\phi$ sur $S$. Remarquons enfin
        que cette construction peut se lisser.
        En faisant de m\^eme sur la branche $S_2$, on r\'ealise le
        rel\`evement de $\phi$ d\'esir\'e. $\Box$\\

\textsc{D. Gayet: Laboratoire Emile Picard, Universit\'e Paul Sabatier,
118 route de Narbonne, 31062 Toulouse Cedex France.} \\
E-mail: gayet@picard.ups-tlse.fr

\end{document}